\documentclass[journal]{IEEEtran}
\usepackage{cite}
\usepackage{amsmath}
\usepackage{amssymb}
\usepackage{amsthm}
\usepackage{arydshln}
\usepackage{mathtools}
\usepackage{mathrsfs}
\usepackage{algorithmic}
\usepackage{graphicx}
\usepackage{epsfig,graphics}
\usepackage{epstopdf}
\usepackage{subfigure}
\usepackage{enumitem}

\usepackage{multirow,ragged2e}
\usepackage{bm}
\usepackage{fixltx2e}
\usepackage[english]{babel}
\usepackage{filecontents}
\theoremstyle{theorem}
\newtheorem{theorem}{Theorem}
\theoremstyle{lemma}

\theoremstyle{definition}

\theoremstyle{assumption}

\theoremstyle{problem}
\newtheorem{problem}{Problem}
\theoremstyle{example}

\newcommand{\fig}[1]{Figure~\ref{#1}}

\newcommand{\sect}[1]{Section~\ref{#1}}
\newcommand{\eq}[1]{Equation~(\ref{#1})}

\newcommand{\kbm}[1]{\bm{\mathfrak{#1}}}
\newcommand{\fbm}[1]{\mathbf{#1}}
\newcommand{\tbm}[1]{\fbm{#1}^\mathsf{T}}
\newcommand{\tfbm}[1]{\bm{#1}^\mathsf{T}}
\newcommand{\ibm}[1]{\fbm{#1}^{-1}}
\newcommand{\itbm}[1]{\fbm{#1}^{-\mathsf{T}}}
\newcommand{\hbm}[1]{\hat{\fbm{#1}}}
\newcommand{\hfbm}[1]{\hat{\bm{#1}}}

\newcommand{\tilfbm}[1]{\tilde{\bm{#1}}}

\newcommand{\ttilfbm}[1]{\tilde{\bm{#1}}^\mathsf{T}}
\newcommand{\dottbm}[1]{\dot{\fbm{#1}}^\mathsf{T}}
\newcommand{\dotbm}[1]{\dot{\fbm{#1}}}

\newcommand{\dothatbm}[1]{\dot{\hat{\fbm{#1}}}}
\newcommand{\dothatfbm}[1]{\dot{\hat{\bm{#1}}}}
\newcommand{\dothattbm}[1]{\dot{\hat{\fbm{#1}}}^\mathsf{T}}

\newcommand{\ddotbm}[1]{\ddot{\fbm{#1}}}
\newcommand{\dddotbm}[1]{\dddot{\fbm{#1}}}
\newcommand{\ddothbm}[1]{\ddot{\hbm{#1}}}
\newcommand{\dddothbm}[1]{\dddot{\hbm{#1}}}

\begin{document}

\title{Bounded Connectivity-Preserving Coordination of Networked Euler-Lagrange Systems}

\author{Yuan~Yang,~\IEEEmembership{Student Member,~IEEE,} Daniela Constantinescu,~\IEEEmembership{Member,~IEEE,} and Yang Shi,~\IEEEmembership{Fellow,~IEEE}\thanks{The authors are with the Department of Mechanical Engineering, University of Victoria, Victoria, BC V8W 2Y2 Canada (e-mail: yangyuan@uvic.ca; danielac@uvic.ca; yshi@uvic.ca).}}

\maketitle

\begin{abstract}
This paper derives sufficient conditions for bounded distributed connectivity-preserving coordination of Euler-Lagrange systems with only position measurements and with system uncertainties, respectively. The paper proposes two strategies that suitably scale conventional gradient-based controls to account for the actuation bounds and to reserve sufficient actuation for damping injection. For output feedback control of networked systems with only position measurements, the paper incorporates a first-order filter to estimate velocities and to inject damping for stability. For networks of uncertain systems, the paper augments conventional linear filter-based adaptive compensation with damping injection to maintain the local connectivity of the network. Analyses based on monotonically decreasing Lyapunov-like functions and Barbalat's lemma lead to sufficient conditions for bounded local connectivity-preserving coordination of Euler-Lagrange networks under the two strategies. The sufficient conditions clarify the interrelationships among the bounded actuations, initial system velocities and initial inter-system distances. Simulation results validate these conditions.
\end{abstract}

\IEEEpeerreviewmaketitle

\section{Introduction}\label{sec:introduction}

Distributed coordination control of multi-agent systems~(MAS) aims to drive each system to the same configuration using only local and $1$-hop state signals~\cite{Ren2005TAC}. Static Proportional~(P) control, for first-order MAS~\cite{Francis2004TAC}, and Proportional-Derivative~(PD) control, for second-order MAS~\cite{Ren2008TAC,Ren2009IJC}, are conventional strategies for coordinating connected MAS. However, because the inter-agent communications of practical MAS-s are distance-dependent, the connectivity assumption may be violated during their coordination~\cite{Zavlanos2011Proceedings}.

Distributed coordination with guaranteed connectivity of first-order MAS-s has been achieved in~\cite{Egerstedt2007TRO,Dimos2007TAC,Dimos2008TRO,Zavlanos2008TRO} through gradient-based controls derived from special unbounded potential functions that quantify the inter-agent coupling energy. Bounded potential functions in~\cite{Amir2010TAC,Dimos2008ICRA} and a dynamic compensation strategy in\cite{Khorasani2013ACC} have accounted for limited actuation and guaranteed connectivity-preserving coordination.

For double-integrator MAS-s, gradient-based controls derived from a bounded potential and the system kinetic energy have preserved connectivity during coordination in~\cite{Su2010SCL}. Adaptation has addressed Lipschitz-like nonlinearities in the system dynamics in~\cite{Su2011Automatica}. Connectivity-preserving leader-follower coordination has been achieved for double-integrator MAS-s in~\cite{Dong2013Automatica,Dong2014TAC,Su2015Automatica,Hu2017TCNS}, and for Euler-Lagrange MAS-s in~\cite{Ren2012SCL,Dong2017IJRNC}. To date, connectivity-preserving coordination of networked Euler-Lagrange systems with actuation constraints has rarely been discussed.

This paper investigates sufficient conditions for guaranteed connectivity-preserving bounded coordination of Euler-Lagrange MAS-s: (i) with only position measurements; and (ii) with uncertain dynamics. A two-agent system illustrates that the actuation bounds and the initial system velocities and inter-system distances determine if the local connectivity of the MAS can be preserved and its coordination achieved. In practice, a coordination strategy designed without considering the actuation constraints may demand larger than available control effort, and thus may not be able to keep the MAS connected to coordinate it. It should be implemented only if it satisfies sufficient conditions for guaranteed connectivity maintenance. Prior to presenting these conditions, the paper:
\begin{enumerate}
\item[1) ]
integrates a first-order filter velocity estimator to achieve output feedback coordination of MAS-s with only position measurements;
\item[2) ]
compensates uncertain dynamics through conventional linear filter-based adaptation, augmented with damping injection to preserve the local network connectivity.
\end{enumerate} 
The paper proposes a general potential function of inter-system distances to characterize the impact of bounded actuation on connectivity preservation. Based on this potential function, the paper derives two sets of sufficient conditions for bounded connectivity-preserving coordination of networked Euler-Lagrange systems, one set for systems that lack velocity sensing and one set for systems with uncertain dynamics.

\section{Problem Statement}\label{sec:problem}
Let a network of non-redundant Euler-Lagrange systems consist of $N$ agents, each with joint space dynamics:
\begin{equation}\label{equ1}
\fbm{M}_{i}(\fbm{q}_{i})\ddotbm{q}_{i}+\fbm{C}_{i}(\fbm{q}_{i},\dotbm{q}_{i})\dotbm{q}_{i}+\fbm{g}_{i}(\fbm{q}_{i})=\tbm{J}_{i}(\fbm{q}_{i})\fbm{f}_{i}\textrm{.}
\end{equation}
In~\eq{equ1}: the subscript $i=1,\cdots,N$ indexes the agent; $\fbm{q}_{i}$, $\dotbm{q}_{i}$ and $\ddotbm{q}_{i}$ are its position, velocity and acceleration; $\fbm{M}_{i}(\fbm{q}_{i})$, $\fbm{C}_{i}(\fbm{q}_{i},\dotbm{q}_{i})$ and $\fbm{J}_{i}(\fbm{q}_{i})$ are its inertia matrix, its matrix of Coriolis and centrifugal effects, and its Jacobian matrix, respectively; $\fbm{g}_{i}(\fbm{q}_{i})$ is the gravitational torque; and $\fbm{f}_{i}$ is the control wrench applied to the agent in task space. 

For Euler-Lagrange systems with only revolute joints, all joint space dynamics in \eq{equ1} have the following property~\cite{Kelly2006Springer}:
\begin{enumerate}[label=P.\arabic*]
\item \label{P1}
The inertia matrices $\fbm{M}_{i}(\fbm{q}_{i})$ are symmetric, positive definite and uniformly bounded by $\fbm{0}\prec\lambda_{i1}\fbm{I}\preceq \fbm{M}_{i}(\fbm{q}_{i})\preceq \lambda_{i2}\fbm{I}\prec\infty$, with $\lambda_{i1}>0, \lambda_{i2}>0$.
\end{enumerate}

This paper addresses the coordination of the end effectors of non-redundant Euler-Lagrange systems away from kinematic singularities, that is, through configurations with full rank Jacobian $\fbm{J}_{i}(\fbm{q}_{i})$. At such configurations, the joint space dynamics are equivalent to the task space dynamics:
\begin{equation}\label{equ2}
\kbm{M}_{i}(\fbm{q}_{i})\ddotbm{x}_{i}+\kbm{C}_{i}(\fbm{q}_{i},\dotbm{q}_{i})\dotbm{x}_{i}+\kbm{g}_{i}(\fbm{q}_{i})=\fbm{f}_{i}\textrm{,}
\end{equation} 
where: $\fbm{x}_{i}$, $\dotbm{x}_{i}$ and $\ddotbm{x}_{i}$ are the end effector configuration, velocity and acceleration in task space; and 
\begin{align*}
\kbm{M}_{i}(\fbm{q}_{i})=&\itbm{J}_{i}(\fbm{q}_{i})\fbm{M}_{i}(\fbm{q}_{i})\ibm{J}_{i}(\fbm{q}_{i})\textrm{,}\\
\kbm{C}_{i}(\fbm{q}_{i},\dotbm{q}_{i})=&\itbm{J}_{i}(\fbm{q}_{i})\fbm{C}_{i}(\fbm{q}_{i},\dotbm{q}_{i})\ibm{J}_{i}(\fbm{q}_{i})\\
&-\itbm{J}_{i}(\fbm{q}_{i})\fbm{M}_{i}(\fbm{q}_{i})\ibm{J}_{i}(\fbm{q}_{i})\dotbm{J}_{i}(\fbm{q}_{i},\dotbm{q}_{i})\ibm{J}_{i}(\fbm{q}_{i})\textrm{,}\\
\kbm{g}_{i}(\fbm{q}_{i})=&\itbm{J}_{i}(\fbm{q}_{i})\fbm{g}_{i}(\fbm{q}_{i})\textrm{.}
\end{align*}
The task-space dynamics~\eqref{equ2} have the following properties~\cite{Khatib1987JRA}:
\begin{enumerate}[label=P.2]
\item \label{P2}
The matrices $\dot{\kbm{M}}_{i}(\fbm{q}_{i})-2\kbm{C}_{i}(\fbm{q}_{i},\dotbm{q}_{i})$ are skew-symmetric.
\end{enumerate}
\begin{enumerate}[label=P.3]
\item \label{P3}
The matrices $\kbm{M}_{i}(\fbm{q}_{i})$ are symmetric, positive definite and uniformly bounded by $\fbm{0}\prec\lambda^{\ast}_{i1}\fbm{I}\preceq\kbm{M}_{i}(\fbm{q}_{i})\preceq\lambda^{\ast}_{i2}\fbm{I}\prec\infty$, with $\lambda^{\ast}_{i1}>0$ and $\lambda^{\ast}_{i2}>0$. 
\end{enumerate}
\begin{enumerate}[label=P.4]
\item \label{P4}
There exist $c_{i}>0$ such that $\|\kbm{C}_{i}(\fbm{q}_{i},\dotbm{q}_{i})\fbm{y}\|\leq c_{i}\|\dotbm{x}_{i}\|\|\fbm{y}\|$, $\forall \fbm{q}_{i}\textrm{,} \dotbm{q}_{i}\textrm{,} \fbm{y}$ such that $\fbm{J}_{i}(\fbm{q}_{i})$ are nonsingular.
\end{enumerate}
\begin{enumerate}[label=P.5]
\item \label{P5}
The networked Euler-Lagrange systems admit a linear parameterization of the form: $\kbm{M}(\fbm{q})\ddotbm{x}_{r}+\kbm{C}(\fbm{q},\dotbm{q})\dotbm{x}_{r}+\kbm{g}(\fbm{q})=\bm{\Phi}(\fbm{q},\dotbm{q},\dotbm{x}_{r},\ddotbm{x}_{r})\bm{\theta}$, where $\bm{\Phi}(\fbm{q},\dotbm{q},\dotbm{x}_{r},\ddotbm{x}_{r})$ is a regressor matrix of known functions and $\underline{\bm{\theta}}\leq\bm{\theta}\leq\overline{\bm{\theta}}$ is a constant vector containing system parameters. 
\end{enumerate} 
For simplicity of notation, matrix and vector dependence on $\fbm{q}_{i}$ and $\dotbm{q}_{i}$ is omitted in the remainder of this paper, for example, $\fbm{C}_{i}(\fbm{q}_{i},\dotbm{q}_{i})$ and $\kbm{g}_{i}(\fbm{q}_{i})$ are indicated as $\fbm{C}_{i}$ and $\kbm{g}_{i}$, respectively.

The communications of the Euler-Lagrange network are described by its communication graph~\cite{Egerstedt2010Princeton}:
\begin{enumerate}[label=D.\arabic*]
\item \label{D1}
The communication graph $\mathcal{G}=\{\mathcal{V},\mathcal{E}\}$ consists of the set of nodes $\mathcal{V}=\{1,\cdots,N\}$ associated with the $N$ systems and the set of edges $\mathcal{E}=\{(j,i)\in\mathcal{V}\times\mathcal{V}|\ j\in\mathcal{N}_{i},\quad \forall i=1,..N\}$ associated with their communications.
\end{enumerate}
In this definition: an edge $(j,i)$ exists if and only if systems $i$ and $j$ exchange state information, that is, they communicate; and $\mathcal{N}_{i}$ is the set of indices of all systems that communicate with, or are adjacent to, system $i$. Graph $\mathcal{G}$ is connected if and only if there exists a path between each pair of systems. A path between systems $i$ and $j$ is a sequence of communication edges of the form $(i,k)$, $(k,h)$, $\cdots$, $(l,m)$ and $(m,j)$. 

The analysis in this paper places the following assumptions on the Euler-Lagrange systems and their communications:
\begin{enumerate}[label=A.\arabic*]
\item \label{A1}
The actuators can more than balance gravity throughout the task space. That is, there exist $\bm{\gamma}_{i}$,~$i=1,...,N$, such that $|\mathfrak{g}^{k}_{i}|\leq \gamma^{k}_{i}<\bar{f}^{k}_{i}\quad \forall k=1,\cdots,n$, where $\bar{\fbm{f}}_{i}$ is the maximum task space actuation wrench of system $i$. Throughout the paper, a $k$ superscript indicates the $k$-th component of a respective vector.
\item \label{A2} No communication edges are built during coordination, $(j,i)\notin\mathcal{E}(0) \Rightarrow (j,i)\notin\mathcal{E}(t)$, $\forall t\geq 0$.
\item \label{A3} The systems have communication radius $r$.
\item \label{A4}
At time $t=0$, the network graph is connected and the distance $d_{ij}=\|\fbm{x}_{ij}\|=\|\fbm{x}_{i}-\fbm{x}_{j}\|$ between the end effectors of the adjacent systems $i$ and $j$ satisfies $d_{ij}\leq r-\epsilon \quad \forall(i,j)\in\mathcal{E}(0)$ for some $\epsilon>0$.
\end{enumerate}

Assumption~\ref{A2} indicates that the communication graph $\mathcal{G}$ of the Euler-Lagrange network is static if the coordination controller preserves connectivity. Assumption~\ref{A3} guarantees that: (i) the distance between the end effectors of any adjacent systems is smaller than the communication radius, $d_{ij}<r$ for all $(i,j)\in \mathcal{E}(t)$ and for any $t\geq 0$; and (ii) the communications are bidirectional, $j\in\mathcal{N}_{i} \Leftrightarrow i\in\mathcal{N}_{j}$, and the graph is undirected. 

For an Euler-Lagrange network with undirected graph, both the weighted adjacency matrix $\fbm{A}_{N\times N}$ and the corresponding weighted Laplacian matrix $\fbm{L}_{N\times N}$ are symmetric, where: $\fbm{A}_{N\times N}=[a_{ij}]$ with $a_{ij}>0$ if $(j,i)\in\mathcal{E}(t)$ and $a_{ij}=0$ otherwise; and $\fbm{L}_{N\times N}=[l_{ij}]$ with $l_{ij}=\sum_{k\in\mathcal{N}_{i}}a_{ik}$ if $j=i$ and $l_{ij}=-a_{ij}$ otherwise. For a network with $M$ pairs of adjacent systems $i$ and $j$, one of the communication edges $(i,j)\in\mathcal{E}(t)$ and $(j,i)\in\mathcal{E}(t)$ is arbitrarily labeled $e_{k}$, $k=1,\cdots,M$, with weight $w(e_{k})=a_{ij}=a_{ji}$. In particular, $e_{k}=(j,i)$ orients the edge from $j$ to $i$, and labels $j$ and $i$ as the tail and the head of $e_{k}$, respectively. After orienting the edges, the incidence matrix is defined by $\fbm{D}_{N\times N}=[d_{hk}]$, with $d_{hk}=1$ if system $h$ is the head of edge $e_{k}$, $d_{hk}=-1$ if system $h$ is the tail of edge $e_{k}$, and $d_{hk}=0$ otherwise. The following lemma gives the relationship between $\fbm{D}$ and $\fbm{L}$:~\cite{Egerstedt2010Princeton}:
\begin{enumerate}[label=L.\arabic*]
\item \label{L1}
The weighted Laplacian matrix $\fbm{L}$ of the undirected connected communication graph $\mathcal{G}$ obeys $\fbm{L}=\fbm{D}\fbm{W}\tbm{D}$, where $\fbm{W}_{M\times M}=\text{diag}\{w(e_{k})\}$.
\end{enumerate}

The problem of connectivity-preserving coordination of networked Euler-Lagrange systems can now be defined.
\begin{problem}\label{pro1}
Given a network of $N$ non-redundant Euler-Lagrange systems that obeys assumptions~\ref{A1}-\ref{A4}, find distributed control laws $\fbm{f}_{i}$ that: 1) maintain all initial communication links and, thus, preserve local connectivity; 2) drive all end effectors to the same task space configuration and, thus, coordinate the systems.
\end{problem}

To explore the feasibility of Problem~\ref{pro1}, consider a simple exemplary network made of two systems with dynamics $\ddot{x}_{i}=f_{i}$, $i=1,2$, communication radius $r>0$, and maximum actuation $\bar{f}_1$ and $\bar{f}_2$. At $t=0$, the systems are connected with $d_{12}(0)=|x_{12}(0)|=r-\epsilon$ for some $0<\epsilon<r$, and move in opposite directions with velocities $\dot{x}_{1}(0)>0$ and $\dot{x}_{2}(0)<0$. To maintain connectivity, the maximum controls should stop the increase of $d_{12}(t)$ before it becomes equal to $r$, or $(|\dot{x}_{1}(0)|+|\dot{x}_{2}(0)|)^{2}<2(r-d_{12}(0))(\bar{f}_{1}+\bar{f}_{2})$. This condition indicates that the feasibility of Problem~\ref{pro1} depends on the bounded actuations. They may be insufficient to maintain the systems in each other's communication radius and, thus, to preserve network connectivity. Sufficient conditions for connectivity-preserving coordination with bounded actuation are, thus, needed to provide a criterion for practical controller implementation.


\section{Main Result}\label{sec:main result}
Let the energy stored in the communication links $(i,j)$ be described by a set of distributed smooth potential functions $\psi(\|\fbm{x}_{ij}\|)$ with the following properties:
\begin{itemize}
\item[1. ]
$\psi(\|\fbm{x}_{ij}\|)$ is positive definite 
for $\|\fbm{x}_{ij}\|\in[0,r]$;
\item[2. ]
$\psi(\|\fbm{x}_{ij}\|)$ is strictly increasing and upper bounded 
on $[0,r]$;
\item[3. ]
$\frac{\partial \psi(\|\fbm{x}_{ij}\|)}{\partial\left(\|\fbm{x}_{ij}\|^{2}\right)}\in [0,\sigma)$ for $\|\fbm{x}_{ij}\|\in[0,r]$;
\item[4. ]
The Hessian of $\psi(\|\fbm{x}_{ij}\|)$ with respect to $\fbm{x}_{ij}$, $\nabla^{2}_{ij}\psi(\|\fbm{x}_{ij}\|)$, is positive definite and upper-bounded by $\nu\fbm{I}$ for $\|\fbm{x}_{ij}\|\in[0,r]$;
\item[5. ]
Given $\|\fbm{x}_{ij}(0)\|\in[0,r-\epsilon]$ and $N>0$ with $0<\epsilon<r$, $N\psi(\|\fbm{x}_{ij}(0)\|)<\psi(r)$ can be guaranteed by tuning function parameters properly. 
\end{itemize}
Here, property 1 converts the coordination problem to the problem of decreasing the potential energy functions to zero. Property 2 makes the preservation of the communication link $(i,j)$ equivalent to guaranteeing $\psi(\|\fbm{x}_{ij}\|)<\psi(r)$. Property 3 leads to a positive and bounded proportional gain of the gradient-based control in the following subsections. Property 4 ensures the existence of a gradient-based velocity control in~\sect{sec: adaptive}. Property 5 facilitates the verification of all communication links through a unifique potential function. A particular function that satisfies the above conditions is $\psi(\|\fbm{x}_{ij}\|)=\frac{\|\fbm{x}_{ij}\|^{2}}{r^{2}-\|\fbm{x}_{ij}\|^{2}+Q}$, with $Q$ a positive constant. 

\subsection{Output Feedback Control}\label{sec: output}
For Euler-Lagrange systems lacking velocity sensors, let a first-order filter provide velocity estimates and the output feedback coordinating controller be designed as:
\begin{equation}\label{equ3}
\begin{aligned}
\hbm{f}_{i}=&-\rho_{i}\sum_{j\in\mathcal{N}_{i}(0)}\nabla_{i}\psi(\|\fbm{x}_{ij}\|)-\kappa_{i}\dothatbm{x}_{i}+\kbm{g}_{i}\textrm{,}\\
\dothatbm{x}_{i}=&-\zeta_{i}\hbm{x}_{i}+\fbm{x}_{i}
\end{aligned}
\end{equation}
where: $\mathcal{N}_{i}(0)$ is the initial set of neighbours of system $i$; $\nabla_{i}\psi(\|\fbm{x}_{ij}\|)=\frac{2\partial\psi(\|\fbm{x}_{ij}\|)}{\partial(\|\fbm{x}_{ij}\|^{2})}(\fbm{x}_{i}-\fbm{x}_{j})$ is the gradient of $\psi(\|\fbm{x}_{ij}\|)$ with respect to $\fbm{x}_{i}$; $\hbm{f}_{i}$ and $\hbm{x}_{i}$ are the computed control force and the estimated velocity of system $i$, respectively; and $\rho_{i}$, $\zeta_{i}$ and $\kappa_{i}$ are positive constants that obey
\begin{equation}\label{equ4}
\begin{aligned}
2|\mathcal{N}_{i}(0)|\rho_{i}\sigma r+\sqrt{2\rho_{i}\kappa_{i}\psi(r)}+\gamma^{k}_{i}\leq \bar{f}^{k}_{i}\quad \forall k=1,..,n \textrm{,}
\end{aligned}
\end{equation} 
with $|\mathcal{N}_{i}(0)|$ the cardinality of $\mathcal{N}_{i}(0)$.

\eq{equ4} guarantees that the bounded actuations can fully implement the scaled gradient-based control $-\rho_{i}\sum_{j\in\mathcal{N}_{i}(0)}\nabla_{i}\psi(\|\fbm{x}_{ij}\|)$ plus damping injection $-\kappa_{i}\dothatbm{x}_{i}$ on $\mathcal{F}_{i}=\left\{(\fbm{x}_{i},\dothatbm{x}_{i})|\ \|\fbm{x}_{ij}\|<r\ \text{for}\ j\in\mathcal{N}_{i}(0)\ \text{and}\ \|\dothatbm{x}_{i}\|<\sqrt{\frac{2\rho_{i}\psi(r)}{\kappa_{i}}}\right\}$. After setting $\hbm{x}_{i}(0)=\frac{\fbm{x}_{i}(0)}{\zeta_{i}}$, assumption~\ref{A4} leads to $(\fbm{x}_{i}(0),\dothatbm{x}_{i}(0))\in \mathcal{F}_{i}$ and further to $\hbm{f}_{i}(0)=\fbm{f}_{i}(0)$, i.e., the initially computed controls can be fully implemented. Then, sufficient conditions for rendering $\mathcal{F}_{i}$ invariant guarantee that $\hbm{f}_{i}(t)=\fbm{f}_{i}(t)\quad \forall t\leq 0$ and the output feedback controller~(\ref{equ3}) achieves connectivity-preserving coordination.


Consider the following potential function for the Euler-Lagrange system \eqref{equ2} under the control \eqref{equ3}:
\begin{equation}\label{equ5}
V=\frac{1}{2}\sum^{N}_{i=1}\left[\frac{1}{\rho_{i}}\left(\dottbm{x}_{i}\kbm{M}_{i}\dotbm{x}_{i}+\kappa_{i}\dothattbm{x}_{i}\dothatbm{x}_{i}\right)+\sum_{j\in\mathcal{N}_{i}(0)}\psi(\|\fbm{x}_{ij}\|)\right]\textrm{.}
\end{equation}
Using $\dotbm{x}_{i}=\fbm{J}_{i}\dotbm{q}_{i}$, property~\ref{P2} and the derivative of the filter dynamics, the derivative of $V$ on $\mathcal{F}=\bigcup_{i=1,..,N}\mathcal{F}_{i}$ is 
\begin{equation}\label{equ6}
\begin{aligned}
\dot{V}=&\frac{1}{2}\sum^{N}_{i=1}\sum_{j\in\mathcal{N}_{i}(0)}\dot{\psi}(\|\fbm{x}_{ij}\|)-\sum^{N}_{i=1}\dottbm{x}_{i}\sum_{j\in\mathcal{N}_{i}(0)}\nabla_{i}\psi(\|\fbm{x}_{ij}\|)\\
&-\sum^{N}_{i=1}\frac{\kappa_{i}}{\rho_{i}}\dottbm{x}_{i}\dothatbm{x}_{i}-\sum^{N}_{i=1}\frac{\kappa_{i}\zeta_{i}}{\rho_{i}}\dothattbm{x}_{i}\dothatbm{x}_{i}+\sum^{N}_{i=1}\frac{\kappa_{i}}{\rho_{i}}\dothatbm{x}_{i}\dotbm{x}_{i}\textrm{,}
\end{aligned}
\end{equation}
Bidirectional communications together with $\psi(\|\fbm{x}_{ij}\|)=\psi(\|\fbm{x}_{ji}\|)$ lead to
\begin{align*}
\frac{1}{2}\sum^{N}_{i=1}\sum_{j\in\mathcal{N}_{i}(0)}\dot{\psi}(\|\fbm{x}_{ij}\|)=\sum^{N}_{i=1}\dottbm{x}_{i}\sum_{j\in\mathcal{N}_{i}(0)}\nabla_{i}\psi(\|\fbm{x}_{ij}\|),
\end{align*}
which, together with \eq{equ6}, leads to 
\begin{equation}\label{equ7}
\dot{V}=-\sum^{N}_{i=1}\frac{\kappa_{i}}{\rho_{i}}\dothattbm{x}_{i}\dothatbm{x}_{i}\leq 0\textrm{.}
\end{equation}

Given that $\|\fbm{x}_{ij}(0)\|\leq r-\epsilon$, $\forall (i,j)\in\mathcal{E}(0)$ by assumption~\ref{A4} and the selection $\hbm{x}_{i}(0)=\frac{\fbm{x}_{i}(0)}{\zeta_{i}}$ makes $\dothatbm{x}_{i}(0)=\fbm{0}$, the conditions
\begin{equation}\label{equ8}
\frac{1}{2}\sum^{N}_{i=1}\left[\frac{\hat{\lambda}_{i20}}{\rho_{i}}\|\dotbm{q}_{i}(0)\|^{2}+\sum_{j\in\mathcal{N}_{i}(0)}\psi(\|\fbm{x}_{ij}(0)\|)\right]<\psi(r)\textrm{,}
\end{equation}
where $\hat{\lambda}_{i20}$ is the maximum eigenvalue of $\fbm{M}_{i}(\fbm{q}_{i}(0))$, imply that $V(t)\leq V(0)<\psi(r)$, i.e., $\psi(\|\fbm{x}_{ij}(t)\|)<\psi(r)$, $\forall (i,j)\in\mathcal{E}(0)$, i.e., local connectivity is preserved. Further, $V(t)<\psi(r)$ implies that $\|\dothatbm{x}_{i}(t)\|<\sqrt{\frac{2\rho_{i}\psi(r)}{\kappa_{i}}}$, $i=1,\cdots,N$, that is, the set $\mathcal{F}$ is invariant and all above statements hold for $t\geq 0$.

Equations \eqref{equ5} and \eqref{equ7} show that $\{\dotbm{x}_{i}\textrm{,}\ \dothatbm{x}_{i}\textrm{,}\ \fbm{x}_{ij}\}\in\mathcal{L}_{\infty}$ and $\dothatbm{x}_{i}\in\mathcal{L}_{2}$ for all $i\in\mathcal{V}\ \text{and}\ (i,j)\in\mathcal{E}(0)$, and together with the derivative of the filter dynamics lead to $\ddothbm{x}_{i}\in\mathcal{L}_{\infty}$ and further to $\dothatbm{x}_{i}\to\fbm{0}$ as $t\to\infty$. After using \eqref{equ1} to show that $\ddotbm{x}_{i}\in\mathcal{L}_{\infty}$, the second-order derivatives of the filter dynamics in~\eqref{equ3} lead to $\dddothbm{x}\in\mathcal{L}_{\infty}$, from which it can be concluded that $\{\ddothbm{x}_{i}\textrm{,}\ \dotbm{x}_{i}\}\to\fbm{0}$ and $\dddotbm{x}_{i}\in\mathcal{L}_{\infty}$, and, by Barbalat's lemma, that $\ddotbm{x}_{i}\to\fbm{0}$. Further, from the system dynamics~\eqref{equ2}, it can be shown that $\sum_{j\in\mathcal{N}_{i}(0)}\nabla_{i}\psi(\|\fbm{x}_{ij}\|)\to\fbm{0}$, $i=1,\cdots,N$. Selecting $\fbm{c}_{k}=[x^{k}_{1}, \cdots, x^{k}_{N}]^\mathsf{T}$, $k=1\cdots,n$, it follows that $\fbm{L}(\fbm{x})\fbm{c}_{k}\to\fbm{0}$, where
\begin{align*}
\fbm{L}(\fbm{x})=[l_{ij}], \quad l_{ij}=\begin{cases}
\sum_{k\in\mathcal{N}_{i}(0)}\frac{\partial\psi(\|\fbm{x}_{ik}\|)}{\partial(\|\fbm{x}_{ik}\|^{2})}\quad &j=i\textrm{,}\\
-\frac{\partial\psi(\|\fbm{x}_{ij}\|)}{\partial(\|\fbm{x}_{ij}\|^{2})}\quad &j\neq i\textrm{.}
\end{cases}
\end{align*}
Lastly, the third property of $\psi(\|\fbm{x}_{ij}\|)$ and lemma~\ref{L1} imply that $\fbm{D}\fbm{c}_{k}\to\fbm{0}$ and $\fbm{x}_{1}\to\fbm{x}_{2}\to\cdots\to\fbm{x}_{N}$, that is, coordination is achieved.

The above proof is summarized in the following theorem:
\begin{theorem}\label{theo1}
Consider the Euler-Lagrange network~\eqref{equ1}, with bounded actuations and satisfying assumptions~\ref{A1}-\ref{A4}. The output feedback control \eqref{equ3} with $\hbm{x}_{i}(0)=\frac{\fbm{x}_{i}(0)}{\zeta_{i}}$ is a connectivity-preserving coordination strategy for the network~\eqref{equ1} if the positive constants $Q$, $\rho_{i}$, $\zeta_{i}$ and $\kappa_{i}$, $i=1,..N$, satisfy conditions \eqref{equ4} and \eqref{equ8}.
\end{theorem} 

\noindent\textbf{Remark. } Condition~\eqref{equ4} guarantees that the designed controls $\hbm{f}_{i}$ can be fully implemented by the limited actuations. Condition~\eqref{equ8} restricts the initial velocities and inter-distances of the Euler-Lagrange systems. Condition~\eqref{equ8} is not practical to check in the absence of velocity measurements. However, for systems initially at rest, condition~\eqref{equ8} can be guaranteed by proper selection of $\psi(\|\fbm{x}_{ij}\|)$. Then, small enough scaling factors $\rho_{i}$ and $\kappa_{i}$ together with assumption~\ref{A1} can limit the designed controls $\hbm{f}_{i}$ to satisfy condition~\eqref{equ4}. Therefore, the solution of the output feedback control problem does always exist for systems initially at rest.


\subsection{Adaptive Control}\label{sec: adaptive}
For systems with uncertain dynamics, the gravity terms cannot be compensated directly. Instead, a linear filter-based adaptive control strategy can be designed as
\begin{equation}\label{equ9}
\begin{aligned}
\hbm{f}_{i}=&\bm{\Phi}_{i}(\fbm{q}_{i},\dotbm{q}_{i},\fbm{e}_{i},\dotbm{e}_{i})\hfbm{\theta}_{i}-\kappa_{i}\fbm{s}_{i}-\mu_{i}\dotbm{x}_{i}\textrm{,}\\
\dothatfbm{\theta}_{i}=&\text{Proj}_{\hfbm{\theta}_{i}}\left(\bm{\omega}_{i}\right)\textrm{,}\\
\bm{\omega}_{i}=&-\beta_{i}\tfbm{\Phi}_{i}(\fbm{q}_{i},\dotbm{q}_{i},\fbm{e}_{i},\dotbm{e}_{i})\fbm{s}_{i}\textrm{,}
\end{aligned}
\end{equation}
where: $i=1,\cdots,N$; $\kappa_{i}$, $\mu_{i}$ and $\beta_{i}$ are positive scaling factors to be determined; $\fbm{s}_{i}=\dotbm{x}_{i}+\alpha\fbm{e}_{i}$ with $\alpha>0$ to be determined and $\fbm{e}_{i}=\sum_{j\in\mathcal{N}_{i}(0)}\nabla_{i}\psi(\|\fbm{x}_{ij}\|)$; and $\bm{\Phi}_{i}(\fbm{q}_{i},\dotbm{q}_{i},\fbm{e}_{i},\dotbm{e}_{i})\hfbm{\theta}_{i}=\hat{\kbm{M}}_{i}(\fbm{q}_{i})[-\alpha\dotbm{e}_{i}]+\hat{\kbm{C}}_{i}(\fbm{q}_{i},\dotbm{q}_{i})[-\alpha\fbm{e}_{i}]+\hat{\kbm{g}}_{i}(\fbm{q}_{i})$. The smooth projection operators $\text{Proj}_{\hfbm{\theta}_{i}}(\bm{\omega}_{i})$ are designed component-wise:
\begin{align*}
\dot{\hat{\theta}}^{k}_{i}=\begin{cases}
\left[1-\upsilon_{lb}(\hat{\theta}^{k}_{i})\right]\omega^{k}_{i}\quad &\underline{\theta}^{k}_{i}\leq\hat{\theta}^{k}_{i}\leq\underline{\theta}^{k}_{i}+\delta\ \&\ \omega^{k}_{i}<0\textrm{,}\\
\left[1-\upsilon_{ub}(\hat{\theta}^{k}_{i})\right]\omega^{k}_{i}\quad &\overline{\theta}^{k}_{i}-\delta\leq\hat{\theta}^{k}_{i}\leq\overline{\theta}^{k}_{i}\ \&\ \omega^{k}_{i}>0\textrm{,}\\
\omega^{k}_{i}\quad &\text{otherwise}\textrm{,}
\end{cases}
\end{align*}
where 
$\upsilon_{lb}(\hat{\theta}^{k}_{i})=\min\left\{1,\frac{\underline{\theta}^{k}_{i}+\delta-\hat{\theta}^{k}_{i}}{\delta}\right\}$ and $\upsilon_{ub}(\hat{\theta}^{k}_{i})=\min\left\{1,\frac{\hat{\theta}^{k}_{i}-\overline{\theta}^{k}_{i}+\delta}{\delta}\right\}$ with $0<\delta<\frac{1}{2}(\overline{\theta}^{k}_{i}-\underline{\theta}^{k}_{i})$, $k=1,\cdots,n$.

From the fourth property of $\psi(\|\fbm{x}_{ij}\|)$, the derivative gains $\nabla^{2}_{ij}\psi(\|\fbm{x}_{ij}\|)$ in $\dotbm{e}_{i}=\sum_{j\in\mathcal{N}_{i}(0)}\nabla^{2}_{ij}\psi(\|\fbm{x}_{ij}\|)(\dotbm{x}_{i}-\dotbm{x}_{j})$ are positive definite and upper-bounded by $\nu\fbm{I}$. From~\cite{Krstic1995}, the projection operators guarantee that $\underline{\bm{\theta}}_{i}\leq\hfbm{\theta}_{i}(t)\leq\overline{\bm{\theta}}_{i}$ by selecting $\underline{\bm{\theta}}_{i}\leq\hfbm{\theta}_{i}(0)\leq\overline{\bm{\theta}}_{i}$, $i=1,\cdots,N$. From~\ref{P1} and~\ref{P4}, if
\begin{equation}\label{equ10}
\begin{aligned}
&\left[(2\nu\lambda^{\ast}_{i2}+c_{i}r\sigma)\alpha|\mathcal{N}_{i}(0)|+\mu_{i}\right]\left[\sqrt{\frac{2\alpha\mu_{i}\psi(r)}{\lambda^{\ast}_{i1}}}+\alpha\sigma r|\mathcal{N}_{i}(0)|\right]\\
&+\gamma^{k}_{i}+\kappa_{i}\sqrt{\frac{2\alpha\mu_{i}\psi(r)}{\lambda^{\ast}_{i1}}}<\bar{f}^{k}_{i}\textrm{,}
\end{aligned}
\end{equation}
hold on $\mathcal{F}_{i}=\Big\{(\fbm{x}_{i},\dotbm{x}_{i})|\ \|\fbm{x}_{ij}\|<r\ \text{for}\ j\in\mathcal{N}_{i}(0)\textrm{,}\ \|\dotbm{x}_{i}\|\leq\sqrt{\frac{2\alpha\mu_{i}\psi(r)}{\lambda^{\ast}_{i1}}}+\alpha\sigma r|\mathcal{N}_{i}(0)|\ \text{and}\ \|\fbm{s}_{i}\|\leq\sqrt{\frac{2\alpha\mu_{i}\psi(r)}{\lambda^{\ast}_{i1}}}\Big\}$, then the adaptive compensation terms $\bm{\Phi}_{i}(\fbm{x}_{i},\dotbm{x}_{i},\fbm{e}_{i},\dotbm{e}_{i})\hfbm{\theta}_{i}$ and the damping injection terms $-\mu_{i}\dotbm{x}_{i}$ can be fully implemented by the limited actuations $\bar{\fbm{f}}_{i}$. When actuators saturate, the remaining actuations can be used for the other terms $-\hat{\kappa}_{i}(t)\fbm{s}_{i}$, where $0<\hat{\kappa}_{i}(t)\leq\kappa_{i}$ are due to actuation limitations. Assuming again that the initial system velocities and inter-distances belong to $\mathcal{F}=\bigcup_{i=1,..,N}\mathcal{F}_{i}$, that is, $\hbm{f}_{i}(0)=\fbm{\Phi}_{i}(\fbm{q}_{i},\dotbm{q}_{i},\fbm{e}_{i},\dotbm{e}_{i})\hfbm{\theta}_{i}-\hat{\kappa}(0)\fbm{s}_{i}-\mu_{i}\dotbm{x}_{i}$, connectivity preservation and the invariance of $\mathcal{F}$ are analyzed using the following Lyapunov-like function:
\begin{equation}\label{equ11}
V=\frac{1}{2}\sum^{N}_{i=1}\left[\frac{1}{\alpha\mu_{i}}\left(\tbm{s}_{i}\kbm{M}_{i}\fbm{s}_{i}+\frac{1}{\beta_{i}}\ttilfbm{\theta}_{i}\tilfbm{\theta}_{i}\right)+\sum_{j\in\mathcal{N}_{i}(0)}\psi(\|\fbm{x}_{ij}\|)\right]\textrm{.}
\end{equation}

After adding $-\bm{\Phi}_{i}(\fbm{x}_{i},\dotbm{x}_{i},\fbm{e}_{i},\dotbm{e}_{i})\bm{\theta}_{i}$ on both sides of~\eqref{equ2} and using~\eqref{equ9}, the system dynamics can be written as
\begin{align*}
&\kbm{M}_{i}(\fbm{q}_{i})\dotbm{s}_{i}+\kbm{C}_{i}(\fbm{q}_{i},\dotbm{q}_{i})\fbm{s}_{i}\\
&=-\bm{\Phi}_{i}(\fbm{q}_{i},\dotbm{q}_{i},\fbm{e}_{i},\dotbm{e}_{i})\tilfbm{\theta}_{i}-\hat{\kappa}_{i}(t)\fbm{s}_{i}-\mu_{i}\dotbm{x}_{i}
\end{align*}
with $i=1,\cdots,N$ and $\tilfbm{\theta}_{i}=\bm{\theta}_{i}-\hfbm{\theta}_{i}$. From~\cite{Krstic1995}, $\ttilfbm{\theta}_{i}[\bm{\omega}_{i}-\dothatfbm{\theta}_{i}]\leq 0$, so the derivative of $V$ on $\mathcal{F}$ is
\begin{align*}
\dot{V}=&\sum^{N}_{i=1}\frac{1}{\alpha\beta_{i}\mu_{i}}\ttilfbm{\theta}_{i}\left(-\beta_{i}\tfbm{\Phi}_{i}(\fbm{x}_{i},\dotbm{x}_{i},\fbm{e}_{i},\dotbm{e}_{i})\fbm{s}_{i}-\dothatfbm{\theta}_{i}\right)-\frac{1}{\alpha}\dottbm{x}_{i}\dotbm{x}_{i}\\
&-\frac{\hat{\kappa}_{i}(t)}{\alpha\mu_{i}}\tbm{s}_{i}\fbm{s}_{i}-\dottbm{x}_{i}\fbm{e}_{i}+\sum^{N}_{i=1}\dottbm{x}_{i}\sum_{j\in\mathcal{N}_{i}(0)}\nabla_{i}\psi(\|\fbm{x}_{ij}\|)\\
\leq&-\sum^{N}_{i=1}\frac{\hat{\kappa}_{i}(t)}{\alpha\mu_{i}}\tbm{s}_{i}\fbm{s}_{i}-\frac{1}{\alpha}\dottbm{x}_{i}\dotbm{x}_{i}\textrm{,}
\end{align*}
which implies that, if conditions~\eqref{equ12} are satisfied
\begin{equation}\label{equ12}
\begin{aligned}
&\frac{1}{2}\sum^{N}_{i=1}\frac{1}{\alpha\mu_{i}}\left(\hat{\lambda}^{\ast}_{i20}\|\fbm{s}_{i}(0)\|^{2}+\frac{1}{\beta_{i}}\|\Delta\bm{\theta}_{i}\|^{2}\right)\\
&+\frac{1}{2}\sum^{N}_{i=1}\sum_{j\in\mathcal{N}_{i}(0)}\psi(\|\fbm{x}_{ij}(0)\|)<\psi(r)\textrm{,}
\end{aligned}
\end{equation}
with $\hat{\lambda}^{\ast}_{i20}$ is the maximum eigenvalue of $\kbm{M}_{i}(\fbm{q}_{i}(0))$ and $\Delta\bm{\theta}_{i}=|\overline{\bm{\theta}}_{i}-\underline{\bm{\theta}}_{i}|$, then $\psi(\|\fbm{x}_{ij}\|)\leq V(t)\leq V(0)<\psi(r)$, i.e., $\mathcal{F}$ is positively invariant and connectivity is guaranteed. 

Further, the derivative of $V$ leads to the conclusions that $\fbm{s}_{i}\in\mathcal{L}_{2}\cap\mathcal{L}_{\infty}$ and $\tilfbm{\theta}_{i}\in\mathcal{L}_{\infty}$, which, together with the system dynamics, lead to $\dotbm{s}_{i}\in\mathcal{L}_{\infty}$ and, thus, to $\fbm{s}_{i}\to\fbm{0}$ as $t\to+\infty$. Using $\fbm{x}=\begin{bmatrix}\tbm{x}_{1},\cdots,\tbm{x}_{N}\end{bmatrix}^\mathsf{T}$ and $\fbm{s}=\begin{bmatrix}\tbm{s}_{1},\cdots,\tbm{s}_{N}\end{bmatrix}^\mathsf{T}$ together with $\dotbm{x}_{i}=-\alpha\fbm{e}_{i}+\fbm{s}_{i}$, the network dynamics become $\dotbm{x}=-\alpha\left[\fbm{L}(\fbm{x})\otimes\fbm{I}_{n}\right]\fbm{x}+\fbm{s}$, where $\fbm{L}(\fbm{x})$ is defined as in~\sect{sec: output}. Then, the analysis procedure in~\cite{Nuno2011TAC} leads to the conclusion that $\fbm{x}_{1}\to\cdots\to\fbm{x}_{N}$ and coordination is achieved.

The above proof is summarized in the following theorem:
\begin{theorem}\label{theo2}
Consider the uncertain Euler-Lagrange network~\eqref{equ1}, with bounded actuations and satisfying assumptions~\ref{A1}-\ref{A4}. The control \eqref{equ9} with $\underline{\bm{\theta}}_{i}\leq\hfbm{\theta}_{i}(0)\leq\overline{\bm{\theta}}_{i}$ is a connectivity-preserving coordination strategy for the network~\eqref{equ1} if the positive constants $Q$, $\alpha$, $\mu_{i}$, $\kappa_{i}$ and $\beta_{i}$, $i=1,..,N$, satisfy conditions \eqref{equ10} and \eqref{equ12}.
\end{theorem}

\noindent\textbf{Remark. } The adaptive terms $\fbm{\Phi}_{i}(\fbm{q}_{i},\dotbm{q}_{i},\fbm{e}_{i},\dotbm{e}_{i})\hfbm{\theta}_{i}$ in~\eqref{equ9} need both joint-space and task-space positions~($\fbm{q}_{i}$ and $\fbm{x}_{i}$) and velocities~($\dotbm{q}_{i}$ and $\dotbm{x}_{i}$). Because system uncertainties thwart the computation of task-space positions and velocities from joint-space positions and velocities through forward kinematics, both the joint-space and task-space positions and velocities should be measured. Further, the upper bounds of $\|\dotbm{x}_{i}\|$ in~\eqref{equ10} are derived from $V(t)<\psi(r)$ and $\fbm{s}_{i}=\dotbm{x}_{i}+\alpha\fbm{e}_{i}$ considering the worst scenarios $\|\dotbm{x}_{i}\|=\|\fbm{s}_{i}\|+\alpha\|\fbm{e}_{i}\|$ because no direct conclusions on $\dotbm{x}_{i}$ bounds can be drawn from $V$. Therefore, conditions~\eqref{equ10} are conservative.

\section{Simulations}\label{sec:simulations}
In this section, two sets of simulations illustrate that bounded connectivity-preserving coordination of Euler-Lagrange networks: (i) cannot be guaranteed if initial system velocities and inter-distances do not meet the sufficient conditions derived in this paper; and (ii) can be achieved if the networked systems initially satisfy the sufficient conditions.

The simulated Euler-Lagrange network comprises $N$ planar two-degrees-of-freedom robots, each with link masses $m_{i1}=m_{i2}=0.5$~kg and lengths $l_{i1}=l_{i2}=1$~m, where the subscript $i=1,\cdots,N$ indexes the robots and subscripts $1$ and $2$ index the robot links. The communication radius is $r=1$~m. The potential function quantifying the energy in each communication link is $\psi(\|\fbm{x}_{ij}\|)=\frac{\|\fbm{x}_{ij}\|^{2}}{r^{2}-\|\fbm{x}_{ij}\|^{2}+Q}$. The control objective is to coordinate all end effectors in task space and preserve all initial communication links in the system with only limited actuations.

\subsection{Connectivity Broken}\label{sec: cb}
This section illustrates that the controllers~\eqref{equ3}, respectively~\eqref{equ9}, cannot maintain the connectivity of a simulated network with $N=2$ robots if their parameters do not satisfy the sufficient conditions~\eqref{equ4} and~\eqref{equ8}, respectively~\eqref{equ10} and~\eqref{equ12}. The robot initial joint space positions and velocities are $\fbm{q}_{1}=[0, \pi/6]^\mathsf{T}$, $\fbm{q}_{2}=[\pi/6, \pi/12]^\mathsf{T}$, $\dotbm{q}_{1}=[-0.5\pi, 0]^\mathsf{T}$ and $\dotbm{q}_{2}=[0.5\pi, 0]^\mathsf{T}$. Selecting $\epsilon=0.2$~m, the robots are initially connected.

\begin{figure}[!hbt]
\centering
\includegraphics[width=\columnwidth , height=5cm]{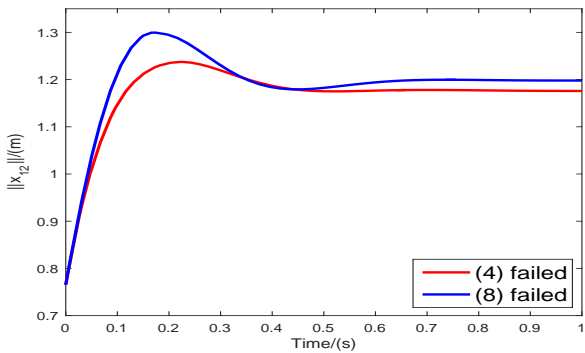}
\caption{Inter-agent distances under output feedback control that fails to satisfy: condition~\eqref{equ4} - red line; (ii) condition~\eqref{equ8} - blue line.}
\label{fig1}
\end{figure}

For output feedback control, the control parameters are chosen as: (a) $Q=0.5$, $\rho_{i}=10$, $\kappa_{i}=100$, $\zeta_{i}=15$ to violate condition~\eqref{equ4} but satisfy condition~\eqref{equ8} given $\bar{f}^{k}_{i}=20$~N; and (b) $Q=5$, $\rho_{i}=10$, $\kappa_{i}=100$, $\zeta_{i}=15$ to violate condition~\eqref{equ8} but satisfy condition~\eqref{equ4} given $\bar{f}^{k}_{i}=40$~N. \fig{fig1} shows that the distance $\|\fbm{x}_{12}\|$ between the two end effectors becomes larger than $r=1$~m, so the communication link breaks in both cases and neither controller coordinates the two end effectors.

\begin{figure}[!hbt]
\centering
\includegraphics[width=\columnwidth , height=5cm]{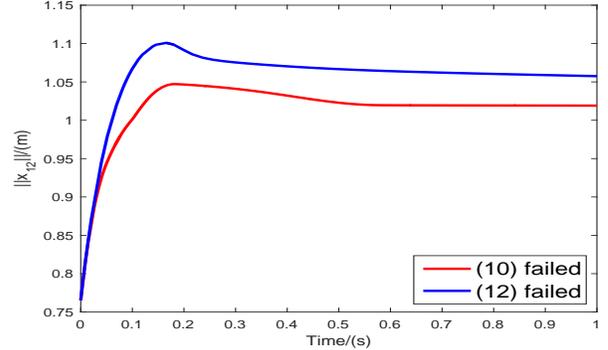}
\caption{Inter-agent distances under adaptive control that fails to satisfy: (i) condition~\eqref{equ10} - red line; (ii) condition~\eqref{equ12} - blue line.}
\label{fig2}
\end{figure}

For adaptive control with actuation bounds $\bar{f}^{k}_{i}=20$~N and with $\|\Delta\bm{\theta}_{i}\|=1.04$, $\beta_{i}=100$ and $\delta_{i}=0.01$, parameters $Q=2$, $\kappa_{i}=10$, $\mu_{i}=400$, $\alpha=10$ satisfy~\eqref{equ12} but not~\eqref{equ10}, and parameters $Q=20$, $\kappa_{i}=10$, $\mu_{i}=5$ and $\alpha=1$ satisfy~\eqref{equ10} but not~\eqref{equ12}. 
 \fig{fig2} illustrates that the distance $\|\fbm{x}_{12}\|$ between the two end effectors increases larger than $r=1$~m. Again, the communication link breaks in both cases and neither controller coordinates the two end effectors.

\subsection{Connectivity Preserved}
This section validates that the controllers~\eqref{equ3}, respectively~\eqref{equ9}, preserve the local connectivity of a simulated network with $N=5$ robots despite bounded actuations if their parameters satisfy the sufficient conditions~\eqref{equ4} and~\eqref{equ8}, respectively~\eqref{equ10} and~\eqref{equ12}. The robots are initially at rest at $\fbm{q}_{1}=[\pi/12, -5\pi/12]^\mathsf{T}$, $\fbm{q}_{2}=[\pi/6, -\pi/3]^\mathsf{T}$, $\fbm{q}_{3}=[\pi/4, -\pi/4]^\mathsf{T}$, $\fbm{q}_{4}=[\pi/3, -\pi/4]^\mathsf{T}$ and $\fbm{q}_{5}=[5\pi/12, -5\pi/12]^\mathsf{T}$. Selecting $\epsilon=0.25$~m guarantees initial connectivity.

\begin{figure}[!hbt]
\centering
\includegraphics[width=\columnwidth , height=6cm]{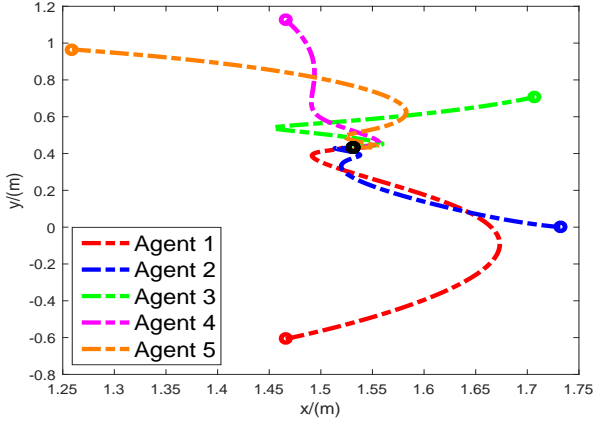}
\caption{End effector task space paths under output feedback control that satisfies the sufficient conditions~\eqref{equ4} and~\eqref{equ8}.}
\label{fig3}
\end{figure}

For output feedback control with $\bar{f}^{k}_{i}=20$~N, the control parameters $Q=1$, $\kappa_{i}=10$, $\zeta_{i}=20$ satisfy the sufficient conditions~\eqref{equ4} and~\eqref{equ8}. The paths of the $5$ end effectors in task space are displayed in~\fig{fig3}. Note in this figure that velocity estimation causes twists and turns during coordination, but the end effectors converge to the same point eventually.  

\begin{figure}[!hbt]
\centering
\includegraphics[width=\columnwidth , height=6cm]{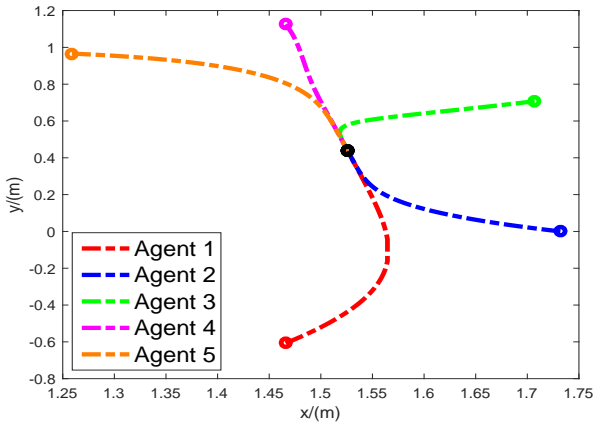}
\caption{End effector task space paths under adaptive control that satisfies the sufficient conditions~\eqref{equ10} and~\eqref{equ12}.}
\label{fig4}
\end{figure}

For adaptive control, conditions~\eqref{equ10} are over conservative due to indirect velocity bounding, and larger actuation bounds are required to identify a controller that satisfies the sufficient conditions~\eqref{equ10} and~\eqref{equ12}. Given $\bar{f}^{k}_{i}=150$~N, for $\|\Delta\bm{\theta}_{i}\|=1.04$ and $\delta_{i}=0.01$, these conditions can be guaranteed by selecting  $Q=0.5$, $\kappa_{i}=100$, $\mu_{i}=10$, $\beta_{i}=500$ and $\alpha=0.001$. \fig{fig4} shows that adaptive control coordinates the $5$ end effectors to the same position as output feedback control, but along smoother paths. Nevertheless, the adaptive compensations make the adaptive coordination significantly slower~($10000$~s)  than the output feedback coordination~($40$~s).

\section{Conclusions}\label{sec:conclusions}
After illustrating that local connectivity preservation is not always achievable in Euler-Lagrange networks with bounded actuations, this paper has explored its sufficient conditions in two cases, output feedback control and adaptive control. For output feedback control, a first-order filter has provided the estimates required for damping injection in systems lacking velocity measurements.
For adaptive control, a linear filter-based adaptive strategy has compensated the uncertain dynamics and has guaranteed the preservation of local connectivity through augmentation with damping injection. The paper has characterized the energy in the communication links through a general class of potential functions that capture the interdependence among bounded actuations, initial system velocities and initial inter-system distances. The analysis of the energy in the communications has led to two theorems summarizing the sufficient conditions for bounded local connectivity-preserving coordination of Euler-Lagrange networks, for output feedback control and for adaptive control, respectively. Because oftentimes networked robots may both lack both velocity measurements and have model uncertainties, future research will investigate the bounded adaptive output feedback coordination of Euler-Lagrange networks with preservation of local connectivity.



\bibliography{bibi}

\begin{thebibliography}{10}
\providecommand{\url}[1]{#1}
\csname url@samestyle\endcsname
\providecommand{\newblock}{\relax}
\providecommand{\bibinfo}[2]{#2}
\providecommand{\BIBentrySTDinterwordspacing}{\spaceskip=0pt\relax}
\providecommand{\BIBentryALTinterwordstretchfactor}{4}
\providecommand{\BIBentryALTinterwordspacing}{\spaceskip=\fontdimen2\font plus
\BIBentryALTinterwordstretchfactor\fontdimen3\font minus
  \fontdimen4\font\relax}
\providecommand{\BIBforeignlanguage}[2]{{%
\expandafter\ifx\csname l@#1\endcsname\relax
\typeout{** WARNING: IEEEtran.bst: No hyphenation pattern has been}%
\typeout{** loaded for the language `#1'. Using the pattern for}%
\typeout{** the default language instead.}%
\else
\language=\csname l@#1\endcsname
\fi
#2}}
\providecommand{\BIBdecl}{\relax}
\BIBdecl

\bibitem{Ren2005TAC}
W.~Ren and R.~W. Beard, ``Consensus seeking in multiagent systems under
  dynamically changing interaction topologies,'' \emph{IEEE Transactions on
  Automatic Control}, vol.~50, no.~5, pp. 655--661, May 2005.

\bibitem{Francis2004TAC}
Z.~Lin, M.~Broucke, and B.~Francis, ``Local control strategies for groups of
  mobile autonomous agents,'' \emph{IEEE Transactions on Automatic Control},
  vol.~49, no.~4, pp. 622--629, April 2004.

\bibitem{Ren2008TAC}
W.~Ren, ``On consensus algorithms for double-integrator dynamics,'' \emph{IEEE
  Transactions on Automatic Control}, vol.~53, no.~6, pp. 1503--1509, July
  2008.

\bibitem{Ren2009IJC}
------, ``Distributed leaderless consensus algorithms for networked
  {E}uler--{L}agrange systems,'' \emph{International Journal of Control},
  vol.~82, no.~11, pp. 2137--2149, 2009.

\bibitem{Zavlanos2011Proceedings}
M.~M. Zavlanos, M.~B. Egerstedt, and G.~J. Pappas, ``Graph-theoretic
  connectivity control of mobile robot networks,'' \emph{Proceedings of the
  IEEE}, vol.~99, no.~9, pp. 1525--1540, Sept 2011.

\bibitem{Egerstedt2007TRO}
M.~Ji and M.~Egerstedt, ``Distributed coordination control of multiagent
  systems while preserving connectedness,'' \emph{IEEE Transactions on
  Robotics}, vol.~23, no.~4, pp. 693--703, Aug 2007.

\bibitem{Dimos2007TAC}
D.~V. Dimarogonas and K.~J. Kyriakopoulos, ``On the rendezvous problem for
  multiple nonholonomic agents,'' \emph{IEEE Transactions on Automatic
  Control}, vol.~52, no.~5, pp. 916--922, May 2007.

\bibitem{Dimos2008TRO}
------, ``Connectedness preserving distributed swarm aggregation for multiple
  kinematic robots,'' \emph{IEEE Transactions on Robotics}, vol.~24, no.~5, pp.
  1213--1223, Oct 2008.

\bibitem{Zavlanos2008TRO}
M.~M. Zavlanos and G.~J. Pappas, ``Distributed connectivity control of mobile
  networks,'' \emph{IEEE Transactions on Robotics}, vol.~24, no.~6, pp.
  1416--1428, Dec 2008.

\bibitem{Amir2010TAC}
A.~Ajorlou, A.~Momeni, and A.~G. Aghdam, ``A class of bounded distributed
  control strategies for connectivity preservation in multi-agent systems,''
  \emph{IEEE Transactions on Automatic Control}, vol.~55, no.~12, pp.
  2828--2833, Dec 2010.

\bibitem{Dimos2008ICRA}
D.~V. Dimarogonas and K.~H. Johansson, ``Decentralized connectivity maintenance
  in mobile networks with bounded inputs,'' in \emph{2008 IEEE International
  Conference on Robotics and Automation}, May 2008, pp. 1507--1512.

\bibitem{Khorasani2013ACC}
I.~Saboori, H.~Nayyeri, and K.~Khorasani, ``A distributed control strategy for
  connectivity preservation of multi-agent systems subject to actuator
  saturation,'' in \emph{2013 American Control Conference}, June 2013, pp.
  4044--4049.

\bibitem{Su2010SCL}
H.~Su, X.~Wang, and G.~Chen, ``Rendezvous of multiple mobile agents with
  preserved network connectivity,'' \emph{Systems \& Control Letters}, vol.~59,
  no.~5, pp. 313 -- 322, 2010.

\bibitem{Su2011Automatica}
H.~Su, G.~Chen, X.~Wang, and Z.~Lin, ``Adaptive second-order consensus of
  networked mobile agents with nonlinear dynamics,'' \emph{Automatica},
  vol.~47, no.~2, pp. 368 -- 375, 2011.

\bibitem{Dong2013Automatica}
Y.~Dong and J.~Huang, ``A leader-following rendezvous problem of double
  integrator multi-agent systems,'' \emph{Automatica}, vol.~49, no.~5, pp. 1386
  -- 1391, 2013.

\bibitem{Dong2014TAC}
------, ``Leader-following connectivity preservation rendezvous of multiple
  double integrator systems based on position measurement only,'' \emph{IEEE
  Transactions on Automatic Control}, vol.~59, no.~9, pp. 2598--2603, Sept
  2014.

\bibitem{Su2015Automatica}
Y.~Su, ``Leader-following rendezvous with connectivity preservation and
  disturbance rejection via internal model approach,'' \emph{Automatica},
  vol.~57, pp. 203 -- 212, 2015.

\bibitem{Hu2017TCNS}
Z.~Feng, C.~Sun, and G.~Hu, ``Robust connectivity preserving rendezvous of
  multirobot systems under unknown dynamics and disturbances,'' \emph{IEEE
  Transactions on Control of Network Systems}, vol.~4, no.~4, pp. 725--735, Dec
  2017.

\bibitem{Ren2012SCL}
Z.~Meng, Z.~Lin, and W.~Ren, ``Leader-follower swarm tracking for networked
  {L}agrange systems,'' \emph{Systems \& Control Letters}, vol.~61, no.~1, pp.
  117 -- 126, 2012.

\bibitem{Dong2017IJRNC}
Y.~Dong and J.~Huang, ``Leader-following consensus with connectivity
  preservation of uncertain {E}uler--{L}agrange multi-agent systems,''
  \emph{International Journal of Robust and Nonlinear Control}, 2017.

\bibitem{Kelly2006Springer}
R.~Kelly, V.~S. Davila, and J.~A.~L. Perez, \emph{Control of robot manipulators
  in joint space}.\hskip 1em plus 0.5em minus 0.4em\relax Springer Science \&
  Business Media, 2006.

\bibitem{Egerstedt2010Princeton}
M.~Mesbahi and M.~Egerstedt, \emph{Graph theoretic methods in multiagent
  networks}.\hskip 1em plus 0.5em minus 0.4em\relax Princeton University Press,
  2010.

\bibitem{Krstic1995}
M.~Krstic, P.~V. Kokotovic, and I.~Kanellakopoulos, \emph{Nonlinear and
  Adaptive Control Design}, 1st~ed.\hskip 1em plus 0.5em minus 0.4em\relax New
  York, NY, USA: John Wiley \& Sons, Inc., 1995.

\bibitem{Nuno2011TAC}
E.~Nu{\~n}o, R.~Ortega, L.~Basa{\~n}ez, and D.~Hill, ``Synchronization of
  networks of nonidentical {E}uler-{L}agrange systems with uncertain parameters
  and communication delays,'' \emph{IEEE Transactions on Automatic Control},
  vol.~56, no.~4, pp. 935--941, April 2011.

\end{thebibliography}
\end{document}